\documentclass[11pt]{amsart}
\usepackage{amsmath,amssymb,amsthm,url}
\usepackage[pdftex]{graphicx}

\DeclareMathOperator{\erfc}{erfc}

\DeclareMathOperator{\majority}{majority}

\newcommand{\beq}{\begin{equation}}
\newcommand{\eeq}{\end{equation}}
\newcommand{\Fix}{\textrm{Fix}}

\renewcommand{\wr}{\mathrm{\,wr\,\,}}

\newtheorem{thm}{Theorem}
\newtheorem{lemma}[thm]{Lemma}
\newtheorem{prop}[thm]{Proposition}
\newtheorem{cor}[thm]{Corollary}

\theoremstyle{definition}

\newtheorem{example}[thm]{Example}

\newtheorem{alg}[thm]{Algorithm}

\begin{document}

\title[Decoding generalised hyperoctahedral groups]{Decoding generalised hyperoctahedral groups and asymptotic analysis of correctible error patterns}

\author{Robert F. Bailey}
\address{Department of Mathematics and Statistics, University of Regina, 3737 Wascana Parkway, Regina, Saskatchewan S4S 0A2, Canada}
\email{robert.bailey@uregina.ca}

\author{Thomas Prellberg}
\address{School of Mathematical Sciences, Queen Mary University of London, Mile End Road, London E1 4NS, United Kingdom}
\email{t.prellberg@qmul.ac.uk}

\begin{abstract} 
We demonstrate a majority-logic decoding algorithm for decoding the
generalised hyperoctahedral group $C_m \wr S_n$ when thought of as
an error-correcting code.  We also find the complexity of this
decoding algorithm and compare it with that of another, more general, algorithm.
Finally, we enumerate the number of error patterns exceeding the
correction capability that can be successfully decoded by this
algorithm, and analyse this asymptotically.\\


\end{abstract}

\maketitle

\section{Introduction} \label{section:intro}

Sets, or groups, of permutations may be used as error-correcting
codes, with permutations in list form as the codewords, and the
usual Hamming distance.  This idea goes back to the 1970s, for
instance to the papers of Blake \cite{Blake74} and Blake, Cohen and
Deza \cite{BlakeCohenDeza79}.  Subsequently, there has been a
resurgence of interest in such codes because of a potential
application to ``powerline communications'', where electrical power
cables are used to transmit data as well as electricity.  For
instance, the 2004 paper by Chu, Colbourn and Dukes
\cite{ChuColbournDukes04} gives a description of this, and some
constructions for suitable codes, while a more general survey can be found
in Huczynska's 2006 paper \cite{Huczynska06}.  More recently, permutation
codes have been applied to ``flash memory'' data storage devices (see the
2010 paper of Tamo and Schwartz \cite{TamoSchwartz10}).

In \cite{ecpg}, the first author gives a decoding algorithm which
works for arbitrary permutation groups when used as codes in this
way.  In \cite{btubb}, he considers certain families
of permutation groups in more detail.  In order to describe the 
case we are interested in in this paper, we need the following 
background material.

Suppose $H$ and $K$ are permutation groups acting on sets $\Gamma$ and $\Delta$ respectively, where $|\Gamma|=m$ and $\Delta=\{1,\ldots,n\}$.  The {\em wreath product} $G=H \wr K$ is constructed as follows.  We consider the action of the Cartesian product $H^n = H\times H \times \cdots \times H$ on $n$ disjoint copies of the set $\Gamma$, labelled by the elements of $\Delta$.  We then form the semidirect product of $H^n$ with $K$, where $K$ acts on $H^n$ according to its action on $\Delta=\{1,\ldots,n\}$; the resulting group is $G=H \wr K$.  Now, we can define an equivalence relation on $\Gamma \times \Delta$ where the equivalence classes are the copies of $\Gamma$; this equivalence relation is preserved by the action of $G$, and so forms a {\em system of imprimitivity} or {\em block system} for $G$.  (See Cameron \cite{CameronPGbook} for more information about permutation groups.)

One family considered in \cite{btubb} were the groups $H \wr S_n$, where $H$ is a {\em regular} permutation group of order $m$.  In this paper, we consider the special case of this where $H$ is a cyclic group of order $m$, so we have $G=C_m \wr S_n$ acting in its imprimitive action on $n$ copies of $\{1,\ldots,m\}$.  We call these groups  {\em generalised hyperoctahedral groups}, as in the case $m=2$ we have the well-known hyperoctahedral group (the automorphism group of the $n$-dimensional hypercube).

In Section \ref{section:algorithm} we give an alternative decoding
algorithm from that given in \cite{ecpg}, that can only be used in
this case, and in Section \ref{section:complexity} we show that this
algorithm is better-performing (in terms of time and space
complexity). Finally, in Section \ref{section:enumeration} we show
that certain patterns of more than $r$ errors can be successfully
decoded by this algorithm, and in Section \ref{section:asymptotics}
we analyse the asymptotic behaviour of this.

Recall that the \emph{minimum distance}, $d$, of a code is the least
value of the Hamming distance over all possible pairs of codewords,
and that the \emph{correction capability} (i.~e.~the number of
errors that can be guaranteed to be corrected), $r$, is given by $r
= \lfloor \frac{d-1}{2} \rfloor$.

\begin{prop} \label{prop:corrcap}
The correction capability of $G=C_m \wr S_n$ is $r=
\lfloor\frac{m-1}{2}\rfloor$.
\end{prop}

\proof
The Hamming distance between two permutations $g$ and $h$ is
precisely $n - |\Fix(gh^{-1})|$, where $\Fix(g)$ is the set of
points fixed by that permutation.  Using the group structure, the
minimum distance of a permutation group $G$ of degree $n$ is
therefore equal to
\[  n - \max_{\substack{ g \in G \\ g \neq id }} \left| \Fix(g) \right| \]
(see \cite{ecpg} for further details).  Now suppose $G=C_m \wr S_n$.
Each copy of $\{1,\ldots,m\}$ forms an imprimitivity block for $G$,
and furthermore if an element of $G$ fixes one point in a block, it
must fix all $m$ points in that block.  Consequently, the maximum
number of fixed points is $(n-1)m$ (it is easy to construct elements
with this many), so the minimum distance is $nm-(n-1)m=m$, and so
the correction capability is $r= \lfloor\frac{m-1}{2}\rfloor$.
\endproof

\section{The decoding algorithm} \label{section:algorithm}

The decoding algorithm we give below makes use of the relatively
straightforward combinatorial structure of the group $G=C_m \wr
S_n$. Since the group permutes $n$ imprimitivity blocks of size $m$,
a permutation in $G$ when written in list form can be divided into
$n$ blocks of length $m$.  The ordering of the blocks gives the action
of $S_n$ on the imprimitivity blocks, and the (cyclic) ordering of the
symbols within a block gives the corresponding element of $C_m$.  Thus each position effectively
holds two pieces of information which are constant throughout that
block: a block label and a cyclic shift.

\begin{example} \label{example:element}
The following permutation is an element of $C_5 \wr S_4$:
\[
[ 7, 8, 9, 10, 6\mid 15, 11, 12, 13, 14\mid 20, 16, 17, 18, 19\mid
5, 1, 2, 3, 4 ].
\]
As can be seen, the permutation splits into four blocks of length
five, containing $1,\ldots,5$, $6,\ldots,10$, etc.
\end{example}

Recall from Proposition \ref{prop:corrcap} that the correction
capability is $r = \lfloor \frac{m-1}{2} \rfloor$.  Consequently, if
we assume there to be a maximum of $r$ errors, there will be a
majority of positions in each block which contain the correct
symbol.  The decoding algorithm uses this fact: the majority of the
block labels and the majority of the cyclic shifts will be correct,
so this allows the reconstruction of the transmitted word.

If a list $L$ has a unique most frequently-occurring element, we
denote it by $\majority(L)$.

\begin{alg} \label{alg:decoding}
Input the received word $w = [w_1|w_2|\cdots|w_n]$, where $w_i =
[w_{i1},\ldots,w_{im}]$.  For each $i$ and $j$, we calculate
$q_{ij}$ and $s_{ij}$ where $w_{ij}=mq_{ij}+s_{ij}$ (where $0\leq q_{ij} \leq n-1$ and $0\leq s_{ij} \leq m-1$), then map
$w_{ij}$ to a pair $(b_{ij},c_{ij})$ as follows:
\begin{eqnarray*}
b_{ij} & := & \left\{ \begin{array}{ll} q_{ij} & \textrm{if $s_{ij} \neq 0$} \\
                                  q_{ij} -1 & \textrm{if $s_{ij} =0$} \end{array}
                                  \right. \\
c_{ij} & := & s_{ij}-j \mod m
\end{eqnarray*}
(where $c_{ij}\in \{0,\ldots,m-1\}$).  Defining
$\widehat{b_i} = \majority [b_{i1},\ldots,b_{im}]$ and $\widehat{c_i} = \majority [c_{i1},\ldots,c_{im}]$
for each $i$, the list
$[\widehat{b_1},\ldots,\widehat{b_n}]$ gives a permutation of
$\{0,1,\ldots,n-1\}$ corresponding to the element of $S_n$ acting on
the blocks, and the list $[\widehat{c_1},\ldots,\widehat{c_n}]$
gives the cyclic shifts within each block.  We can then reconstruct
the original permutation $g := [g_1|g_2|\cdots|g_n]$, where $g_i =
[g_{i1},\ldots,g_{im}]$, $g_{ij} = m\widehat{b_i} + t_{ij}$ and
$t_{ij}= j + \widehat{c_i} \mod m$.  (Note that we assume $t_{ij}\in \{1,\ldots,m\}$\footnote{This is done to reconcile two conventions, namely that permutations are of the set $\{1,\ldots,m\}$ while modular arithmetic is performed on the set $\{0,\ldots,m-1\}$.}.)
\end{alg}

\begin{example} \label{example:algexample}
Suppose we transmit the following element \mbox{$g \in C_5 \wr S_4$}:
\[
[ 7, 8, 9, 10, 6\mid 15, 11, 12, 13, 14\mid 20, 16, 17, 18,
19\mid 5, 1, 2, 3, 4 ].
\]
Then suppose we receive the following word $w$:
\[
[17, 1, 9, 10, 6\mid 15, 11, 12, 13, 14\mid 20, 16, 17, 18,
19\mid 5, 1, 2, 3, 4].
\]
This clearly has errors in positions 1 and 2.  Having split this
into four blocks of length five, we obtain the data shown in Table 
\ref{table:data}.

\begin{table}[hbtp]
\centering
\[
\begin{array}{|cc|c|c|}

\hline i & j & w_{ij} & (b_{ij},c_{ij}) \\ \hline

0 & 1 & 17 & ( 3, 1 ) \\
0 & 2 & 1  & ( 0, 4 ) \\
0 & 3 & 9  & ( 1, 1 ) \\
0 & 4 & 10 & ( 1, 1 ) \\
0 & 5 & 6  & ( 1, 1 ) \\ \hline

1 & 1 & 15 & ( 2, 4 ) \\
1 & 2 & 11 & ( 2, 4 ) \\
1 & 3 & 12 & ( 2, 4 ) \\
1 & 4 & 13 & ( 2, 4 ) \\
1 & 5 & 14 & ( 2, 4 ) \\ \hline
\end{array}
\qquad 
\begin{array}{|cc|c|c|}

\hline i & j & w_{ij} & (b_{ij},c_{ij}) \\ \hline
2 & 1 & 20 & ( 3, 4 ) \\
2 & 2 & 16 & ( 3, 4 ) \\
2 & 3 & 17 & ( 3, 4 ) \\
2 & 4 & 18 & ( 3, 4 ) \\
2 & 5 & 19 & ( 3, 4 ) \\ \hline

3 & 1 & 5  & ( 0, 4 ) \\
3 & 2 & 1  & ( 0, 4 ) \\
3 & 3 & 2  & ( 0, 4 ) \\
3 & 4 & 3  & ( 0, 4 ) \\
3 & 5 & 4  & ( 0, 4 ) \\ \hline

\end{array}
\]
\caption{ \label{table:data} Data obtained during decoding in Example \ref{example:algexample}}
\end{table}

Taking the ``majority'' elements, we find the block permutation
$\beta = [ 1, 2, 3, 0 ]$ and cyclic shifts of $[ 1, 4, 4, 4 ]$.  We
have the information needed to reconstruct the transmitted
permutation: for instance with $i=0$ and $j=1$, we have
$\widehat{b_0}=1$, $\widehat{c_0}=1$, $t_{01}= 1+1 \mod 5 = 2$ and
so $g_{01}=1\times 5 + 2 =7$.  Performing these calculations for
each $i$ and $j$, we can recover the transmitted permutation:
\[
[ 7, 8, 9, 10, 6\mid 15, 11, 12, 13, 14\mid 20, 16, 17, 18, 19\mid
5, 1, 2, 3, 4 ].
\]
\end{example}

We conclude this section by mentioning that the first author has 
implemented this algorithm in the computer algebra system {\sf GAP} \cite{GAP}.

\section{Complexity} \label{section:complexity}
We recall that Algorithm~\ref{alg:decoding} has three parts:
calculating the numbers $(b_{ij},c_{ij})$, which involves integer
arithmetic; determining the most frequently-occurring elements
$\widehat{b_i}$ and $\widehat{c_i}$; and reconstructing the decoding
permutation $g$, which involves more integer arithmetic.  In order
to determine the complexity of this algorithm, there are some
assumptions we need to make first.
\begin{itemize}
\item integer arithmetic can be done via a look-up table, in constant time;
\item comparing the sizes of two integers can be done in constant time;
\item finding position $i$ in a list of length $k$ takes $\mathrm{O}(\log k)$ time.
\end{itemize}
However, the second step is more complicated and requires the
following lemma.

\begin{lemma} \label{lemma:majority}
Let $L$ be a list of length $m$ with symbols chosen from
$S=\{1,\ldots,k\}$.  Suppose $L$ has a unique most frequently
occurring element, $x \in S$.  Then the time taken to determine $x$
is $\mathrm{O}(k + m \log k)$.
\end{lemma}

\proof  We begin by producing an auxiliary list $K$ of length $k$, initially set to $[0,0,\ldots,0]$.  This takes $k$ units of time.  We then work through each of the positions of $L$: in each position, we do as follows:
\begin{itemize}
\item read the symbol, $i$, (taking one unit of time);
\item find position $i$ in $K$ (taking $\mathrm{O}(\log k)$ time);
\item increment that entry by 1 (taking one unit of time).
\end{itemize}
This turns $K$ into a list of the frequencies of each symbol in $S$
in the list $L$.  Doing this for each of the $m$ entries of $L$
requires a total of $\mathrm{O}(m \log k)$ time units.  We then work
through $K$ to find the position of the maximum element.  This will
require $\mathrm{O}(k)$ comparisons.  Combining this, we have
$\mathrm{O}(k) + \mathrm{O}(m \log k) + \mathrm{O}(k) = \mathrm{O}(k
+ m \log k)$ as required.  \endproof

We observe that the method described above is not necessarily the
best possible; other methods may be faster, and which method is the
best may depend on factors such as the relative sizes of $m$ and
$k$.  However, we now use it to determine the time complexity of
Algorithm \ref{alg:decoding}.

\begin{thm} \label{thm:time}
The time required to perform Algorithm \ref{alg:decoding} is
$\mathrm{O}(mn \log m)$ (if $m \geq n$) or $\mathrm{O}(n^2 + mn\log
n)$ (if $m \leq n$).
\end{thm}

\proof
The first stage is the calculation of the numbers $(b_{ij},c_{ij})$.
There are $mn$ such calculations to perform, and we have assumed
that each takes a constant amount of time, requiring a total of
$\mathrm{O}(mn)$ time units.

The next stage is, in each block $i$, to determine the most
frequently occurring block label $\widehat{b_i}$ and most frequently
occurring cyclic shift $\widehat{c_i}$.  This involves determining
the most frequently occurring element in a list of length $m$ with
symbols chosen from a set of size $n$ (for the block labels) and
from a list of length $m$ with symbols chosen from a set of size $m$
(for the cyclic shifts).  By Lemma \ref{lemma:majority} above, the
first of these will take $\mathrm{O}(n+m \log n)$ time units, the
second $\mathrm{O}(m + m \log m)$.  As this has to be done in each
of $n$ blocks, this gives a total of $\mathrm{O}(n^2 + mn \log n +
mn + mn \log m)$.

The final stage is the reconstruction of $g$, which requires $m$
integer arithmetic operations in each of the $n$ blocks.  As we have
assumed that integer arithmetic takes constant time, this requires a
total of $\mathrm{O}(mn)$ time.  So the total time required is
$\mathrm{O}(mn) + \mathrm{O}(n^2 + mn \log n + mn + mn \log m) +
\mathrm{O}(mn)$.  If $m \geq n$, this reduces to $\mathrm{O}(mn + mn
\log m) = \mathrm{O}(mn \log mn)$, while if $m \leq n$ it reduces to
$\mathrm{O}(n^2 + mn\log n)$, as required.
\endproof

We should also consider the space complexity of Algorithm
\ref{alg:decoding}.  This time, we require a look-up table for our
integer arithmetic, and there are also items that have to be stored
whilst the algorithm is being performed.

\begin{prop} \label{prop:space}
The amount of storage space required by the decoding algorithm is
$\mathrm{O}(mn^2)$, and the space required to perform the algorithm
is $\mathrm{O}(mn)$.
\end{prop}

\proof We need to store a look-up table, where for each of the $mn$ symbols, for $n$ possible divisors we record a quotient/remainder pair.  This requires a total of $2mn^2 = \mathrm{O}(mn^2)$ storage units.  To perform the algorithm, we need to store $mn$ quotient/remainder pairs, then the two auxiliary lists (one of length $m$ and one of length $n$) to find their most frequently-occurring element, and need $mn$ units to store the reconstructed group element.  This gives a total of $2mn + m + n + mn = \mathrm{O}(mn)$ units.  \endproof

In \cite{ecpg}, a more general decoding algorithm was given, which works for arbitrary permutation groups; also in \cite{ecpg} its complexity was analysed in a similar fashion to the above.  In the case where the group is the generalised hyperoctahedral group $C_m \wr S_n$, bounds on the time and space complexity for the more general algorithm are given by $\mathrm{O}(m^2 n^2)$ and $O(m^3 n^3)$ respectively.  Thus, from the point of view of a worst-case analysis, Theorem \ref{thm:time} and Proposition \ref{prop:space} suggest that Algorithm \ref{alg:decoding} is an improvement.  (Both algorithms require $\mathrm{O}(mn)$ space to perform the algorithm.)

\section{Enumerating correctible error patterns} \label{section:enumeration}

Suppose we have transmitted a permutation $g\in C_m\wr S_n$ and obtained the received word~$w$, which contains errors.  The {\em error pattern} of~$w$ is the subset of the positions $\{1,\ldots,mn\}$ where the errors are situated.  Formally, a {\em $k$-error pattern} is a subset of $\{1,\ldots,mn\}$ of size $k$.

We observe that Algorithm~\ref{alg:decoding} will successfully decode~$w$ if there are a majority of correct elements \emph{in each block}.  Consequently, there will be error patterns of size at most~$nr$ that can successfully be corrected (where $r = \lfloor \frac{m-1}{2} \rfloor$), regardless of what the erroneous symbols are.  We call an error pattern {\em correctible} if it contains no more than~$r$ errors in each block.  In this section, we investigate how many such patterns there are.

Before we do so, we remark that for a given transmitted permutation~$g$ there are received words whose error patterns are not correctible, but which still can be decoded by Algorithm~\ref{alg:decoding}, depending on the nature of the erroneous symbols.  For instance, consider Example~\ref{example:algexample}, but suppose the received word $w$ begins $[7,8,6,15,1\mid \ldots]$.  Three positions in that block contain errors (so the pattern is not correctible), yet Algorithm~\ref{alg:decoding} would determine the correct block label and cyclic shift.  On the other hand, if the received word $w$ begins $[7,8,1,1,1\mid\ldots]$, the error pattern is the same, but Algorithm~\ref{alg:decoding} would fail.  In the remainder of the paper, we are only concerned with correctible error patterns.

For positive integers $k$, $n$ and $r$, define $\mathcal{P}_{n,r}(k)$ to be 
the set of all partitions of the integer $k$ into at most $n$ parts, and
where each part has size at most $r$.  For $\pi \in \mathcal{P}_{n,r}(k)$, 
we denote the number of parts of size $i$ by $f_i(\pi)$ (so that 
$\sum f_i(\pi) \leq n$).  We also define a quantity $c_i(\pi)$ to be
\[ c_i(\pi) = \sum_{j=1}^{i-1} f_j(\pi)\]
for $i\geq2$, with $c_1(\pi)=0$.  That is, $c_i(\pi)$ is the number of parts 
in $\pi$ of size strictly less than $i$.

\begin{prop} \label{prop:counting}
For a word in $C_m \wr S_n$, and for $k \leq nr$, the number of $k$-error patterns which are correctible is given by
\[ E_{n,m,r}(k) = \sum_{\pi \in \mathcal{P}_{n,r}(k)} \prod_{i=1}^r {n-c_i(\pi) \choose f_i(\pi)} {m \choose i}^{f_i(\pi)}.\]
\end{prop}

\proof For a $k$-error pattern to be correctible, the errors can be spread across up to $n$ blocks, as long as there are no more than $r$ errors in each block.  So for a given partition $\pi \in \mathcal{P}_{n,r}(k)$, for each part of $\pi$ we have to choose (i) which block contains that many errors and (ii) where in that block they lie.  Working through $i$ in increasing order, for each $i$ there are $n-c_i(\pi)$ blocks remaining, of which we choose $f_i(\pi)$ (corresponding to the $f_i(\pi)$ parts of size $i$).  Then in each of the $f_i(\pi)$ blocks we have chosen, we choose $i$ error positions from the $m$ available. \endproof

\noindent Note that if $k > nr$, the set $\mathcal{P}_{n,r}(k)$ is
empty, so $E_{n,m,r}(k)=0$.

While this is a tidy combinatorial expression for the desired
quantity $E_{n,m,r}(k)$, its behaviour cannot easily be seen,
especially as we wish to compare it with the total number of
$k$-error patterns~${ mn \choose k }$.  A first step would be
to find a recurrence relation.

\begin{lemma} \label{prop:recurrence}
The numbers $E_{n,m,r}(k)$ satisfy the recurrence relation
\[ E_{n,m,r}(k) = \sum_{l=0}^r {m \choose l} E_{n-1,m,r}(k-l).\]
\end{lemma}

\proof Suppose there are $l$ errors in the $n^\textrm{th}$ block; there are ${m \choose l}$ ways of arranging these.  Then there are $k-l$ errors in the remaining $n-1$ blocks, so there are $E_{n-1,m,r}(k-l)$ ways of arranging these.  Summing over all possible values of $l \leq r$, we obtain the required relation. \endproof

This recurrence relation assists us in studying the generating function for $E_{n,m,r}(k)$.  Let $\mathcal{E}_{n,m,r}(x)$ denote this function, that is
\[ \mathcal{E}_{n,m,r}(x) = \sum_{k \geq 0} E_{n,m,r}(k)\ x^k,\]
and define \[ F_{m,r}(x)=\sum_{l=0}^r {m \choose l} x^l.\]

\begin{prop} \label{prop:genfuncrecurrence}
The generating function $\mathcal{E}_{n,m,r}(x)$ can be rewritten as
$$\mathcal{E}_{n,m,r}(x) = (F_{m,r}(x))^n\;.$$
\end{prop}

\proof
Applying Lemma \ref{prop:recurrence} and re-summing, we obtain
\begin{eqnarray*}
\mathcal{E}_{n,m,r}(x) & = & \sum_{k \geq 0} E_{n,m,r}(k) x^k \\
                 & = & \sum_{k \geq 0} \sum_{l=0}^r {m \choose l} E_{n-1,m,r}(k-l) x^{k-l} x^l \\
                 & = & \left( \sum_{l=0}^r {m \choose l} x^l \right) \mathcal{E}_{n-1,m,r}(x).
\end{eqnarray*}
By iterating this, and with the observation that $\mathcal{E}_{0,m,r}(x)=1$, we have
\[ \mathcal{E}_{n,m,r}(x) = \left( \sum_{l=0}^r {m \choose l} x^l \right)^n \]
as required.
\endproof

The probability that a $k$-error pattern is correctible is then given by
\[
p_{n,m,r}(k)=
E_{n,m,r}(k)/\textstyle{nm\choose k}\;,
\]
as the number of all possible $k$-error patterns is~${nm\choose k}$, of which $E_{n,m,r}(k)$ are correctible.

From the point of view of applications it is perhaps more useful to consider the
probability $P_{n,m,r}(p)$ that a received word has a correctible error pattern,
under the assumption that an individual error occurs with probability $p$. In other words, we consider a probabilistic model in which the number $k$ of errors is binomially distributed, $k\sim B(mn,p)$, so that the expected number of errors is given by $pmn$.

\begin{prop} \label{prop:expectederror}
For a word in $C_m \wr S_n$, if an individual error occurs
with probability $p$ then the probability that its error pattern is correctible is given by
\[
P_{n,m,r}(p)=(1-p)^{mn}\mathcal{E}_{n,m,r}\left(p/(1-p)\right)\;.
\]
\end{prop}

\proof
An error pattern is correctible if there are at most $r$ errors in each block. The probability that
exactly $l$ errors occur in a block of size $m$ is given by ${m\choose l}p^l(1-p)^{m-l}$, so
that the probability that there are at most $r$ errors in each of $n$ blocks is given by
\[
P_{n,m,r}(p)=\left(\sum_{l=0}^r {m \choose l}p^l(1-p)^{m-l}\right)^n\;,
\]
which gives the desired expression.
\endproof

Of course these two probabilistic models are not equivalent.

\section{Asymptotic analysis} \label{section:asymptotics}

\subsection{Asymptotics of $p_{n,m,r}(k)$ for $C_m\wr S_n$ as $n\rightarrow\infty$}

Let us focus on the asymptotics of $p_{n,m,r}(k)$ when the error frequency $k/mn$ is fixed.
Here, we will discuss the case where $n\rightarrow\infty$ as $m$ is fixed.

We first give an expression for $E_{n,m,r}(k)$ which is amenable to asymptotic treatment.

\begin{lemma} \label{lemma:cauchyformula}
\[
E_{n,m,r}(k)=\frac1{2\pi i}\oint\left(F_{m,r}(z)\right)^n\frac{dz}{z^{k+1}}\;,
\]
where the contour of integration is a counterclockwise circle about the origin.
\end{lemma}

\proof
This follows directly from Proposition \ref{prop:genfuncrecurrence} and the Cauchy Integral Formula.
\endproof

An asymptotic analysis of the integral in Lemma
\ref{lemma:cauchyformula} is obtained from a saddle-point
approximation.  (See Flajolet and Sedgewick \cite{Flajolet} for
background material on this technique.)  For $n$ and $k$ large, the
behaviour of the integral is determined by the exponential of
\[
n\log F_{m,r}(z)-k\log z\;.
\]
There is a unique positive saddle $\zeta$ given by
\[
\label{saddle} 0=\frac d{d\zeta}\left[\log F_{m,r}(\zeta)-\frac
kn\log\zeta\right]\;,
\]
and the asymptotics are obtained by approximating the integrand
around the saddle by a Gaussian. This is the content of Theorem
VIII.8 of \cite{Flajolet}, which we use to obtain the following
result.

\begin{prop} \label{prop:asymptotics1}
Let $\lambda=k/n$ be a fixed positive number with \mbox{$0<\lambda<r$}, let $\zeta$ be the
unique positive root of the equation
\[\zeta\frac{F_{m,r}'(\zeta)}{F_{m,r}(\zeta)}=\lambda\;,\]
and let
\[
\xi=\frac{d^2}{d\zeta^2}\left[\log F_{m,r}(\zeta)-\frac kn\log\zeta\right]\;.
\]
Then, with $k=\lambda n$ an integer, one has, as $n\rightarrow\infty$,
\[
E_{n,m,r}(k)=\frac{F_{m,r}(\zeta)^n}{\zeta^{k+1}\sqrt{2\pi n\xi}}(1+o(1))\;.
\]
In addition, a full expansion in descending powers of $n$ exists.
These estimates hold uniformly for $\lambda$ in any compact interval of $[0,r]$.
\end{prop}

\proof
One easily checks that the conditions of Theorem VIII.8 in \cite{Flajolet} are satisfied.
\endproof

Fixing $m$ and the fraction of errors $k/mn=\lambda/m$, this allows us to control the 
asymptotics of $p_{n,m,r}(k)$ for large $n$.

\begin{figure}[htb]  
\vspace*{-2cm}
\begin{center}\includegraphics[angle=0,width=0.7\textwidth]{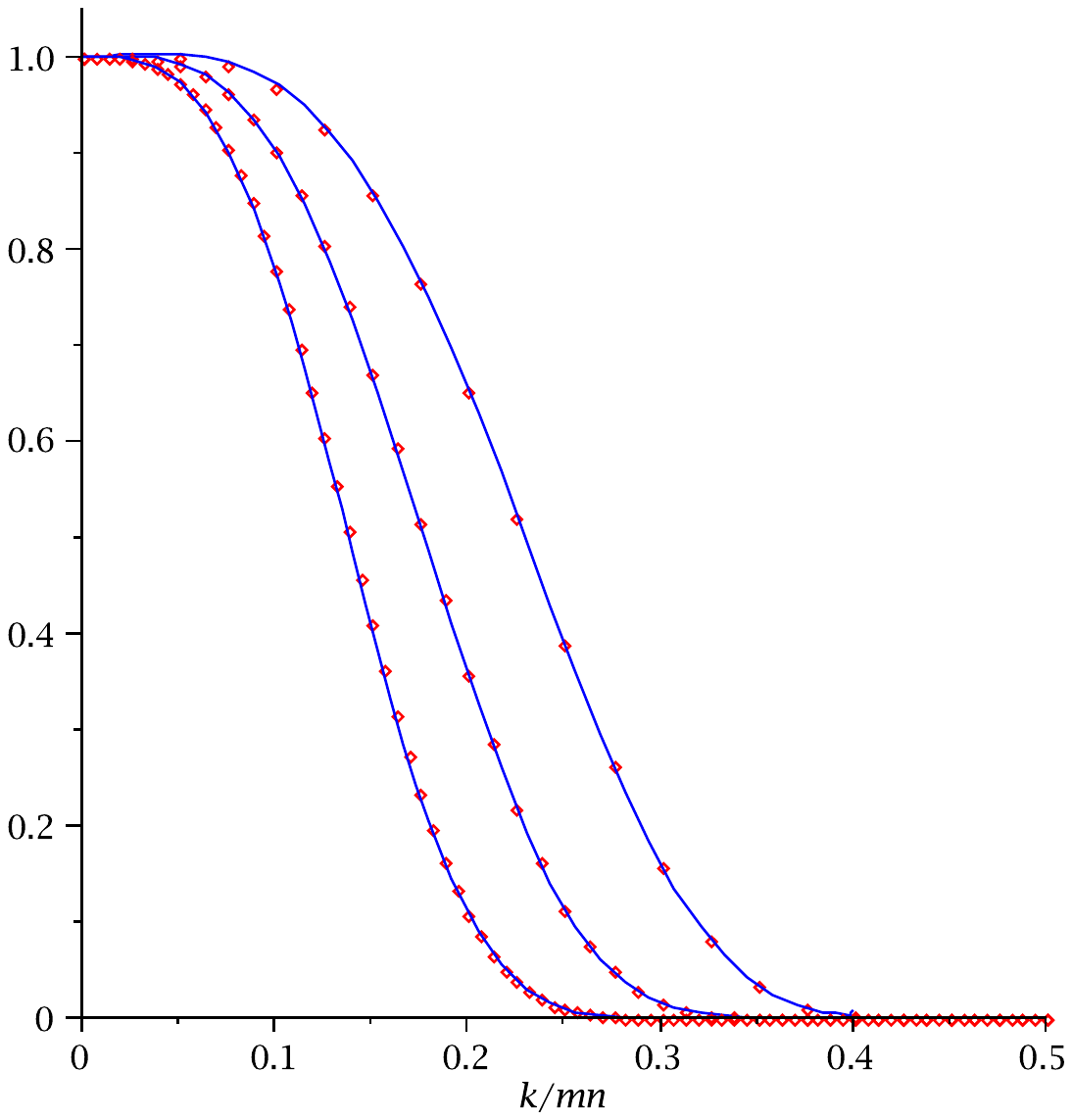}\end{center}
\vspace*{-2cm}
\caption{The probability $p_{n,m,r}(k)$ that a $k$-error pattern is correctible versus the error frequency $k/mn$ for $C_m\wr S_n$, when $m=5$ (and thus $r=\lfloor\frac{5-1}2\rfloor=2$).
Shown is a comparison of the exact values (shown as points) and the asymptotic result (shown as curves) from Proposition \ref{prop:asymptotics1} for $n=8$, $16$, and $32$ (from right to left).} \label{fig52}
\end{figure}

\begin{example} \label{example:asy1example} 
Figure~\ref{fig52} shows the probability $p_{n,m,r}(k)$ that a $k$-error pattern is correctible for $C_m\wr S_n$, for three different values
of $n$ and when $m=5$ (and thus $r=\lfloor\frac{5-1}2\rfloor=2$). 
To reduce the error in the asymptotic formula for small values of $k$, we 
replace both the numerator $E_{n,m,r}(k)$ and denominator ${nm\choose k}$ of $p_{n,m,r}(k)$
by the respective leading terms of the asymptotic expansion given by Proposition~\ref{prop:asymptotics1} (for the denominator we use $r=m$). One
expects heuristically that first-order corrections to the leading asymptotics will largely
cancel each other. Numerically this seems to be confirmed, as even for moderate values 
of~$n$ the agreement between the asymptotic result and the exact values is remarkably good.
\end{example}

\subsection{Asymptotics of $P_{n,m,r}(p)$ for $C_m\wr S_n$ as $m\rightarrow\infty$}

Let us now consider the asymptotics of $P_{n,m,r}(p)$, given that the probability $p$ of a single error is fixed. Here, we will discuss the case of $m\rightarrow\infty$ as $n$ is fixed.

To deal with the truncated binomial sum $F_{m,r}(x)$, we will use the following integral formulation.

\begin{lemma} \label{lemma:Fintegral}
\[ F_{m,r}(x)=\frac1{2\pi i}\oint{(1+sx)^m}\frac{ds}{s^{r+1}(1-s)}\;,\]
where the contour of integration is a clockwise circle about the origin of radius
less than one.
\end{lemma}

\proof
Expand the integrand as a power series in $s$ (which is absolutely convergent for $|s|<1$)
and integrate term-by-term.
\endproof

As in the previous subsection, an asymptotic analysis of the
integral in Lemma \ref{lemma:Fintegral} is obtained from a
saddle-point approximation.  For $m$ and $r$ large, the behaviour of
the integral is determined by the exponential of
\[ m\log(1+sx)-r\log s\;.\]
There is a unique saddle $\sigma$ given by
\[
0=\frac d{d\sigma}\left[\log(1+\sigma x)-\frac rm\log\sigma\right]\;,
\]
namely,
\[
\sigma=\frac1x\frac r{m-r}\;.
\]

Importantly, the saddle collides with the amplitude critical point
$s=1$ when $x=r/(m-r)$. This changes the asymptotic behaviour
significantly, and we need a {\em uniform} asymptotic expansion to
take this into consideration. The standard procedure here is to
re-parameterise the contour by a quadratic, i.e.,
\[
\log(1+sx)-\frac rm\log s=-\frac{t^2}2-\gamma t+\delta\;,
\]
where $\gamma$ and $\delta$ are determined by matching the location of the saddle point
$s=\sigma$ with
$t=-\gamma$, and the location of the critical point $s=1$ with $t=0$.

In the following result, $\erfc$ denotes the {\em complementary
error function}, which is defined as
\[ \erfc(x) = \frac{2}{\sqrt{\pi}} \int_x^{\infty} e^{-t^2}\,dt.\]
We shall also need
\[
\rho(\beta,x)=\sqrt{\log(1+x)-\beta\log x-h(\beta)}
\]
where
\[
h(\beta)=-\beta\log\beta-(1-\beta)\log(1-\beta)\;,
\]
and
\[
A(\beta,x)=\frac1{\sqrt{\beta(1-\beta)}\left(1-\frac\beta{x(1-\beta)}\right)}-\frac1{\sqrt2\rho(\beta,x)}\;.
\]
Note that the radicand in $\rho(\beta,x)$ has a quadratic zero at \mbox{$x=\beta/(1-\beta)$}, and
that the correct sign has to be chosen to make $\rho(\beta,x)$ an analytic function
near that point.

\begin{prop} \label{prop:asymptotics2}
Let $\beta=r/m$ be a fixed positive number with \mbox{$0<\beta<1$.}
Then, with $r=\beta m$ an integer, we have, as $m\rightarrow\infty$,
\[
F_{m,r}(x)=%&
(1+x)^m\left[\frac12\erfc(\sqrt{m}\rho(\beta,x))\right.%\\
+\left.\frac{A(\beta,x)}{\sqrt{2m\pi}}e^{-m\rho(\beta,x)^2}\right](1+o(1))\;.
\]
In addition, a full expansion in descending powers of $m$ exists.
These estimates hold uniformly for $\beta$ in any compact interval of $[0,1]$ and $x$
in any compact domain of $\mathbb C\setminus \mathbb R^-$.
\end{prop}

\proof
This result follows from equation (9.4.22) in Section 9.4 of
Bleistein and Handelsman \cite{Bleistein} (with $r=0$).
\endproof

As an alternative to using Lemma \ref{lemma:Fintegral}, one could
have written $F_{m,r}(x)$ in terms of an incomplete Beta function and used
results of Temme \cite{Temme}.

The main result now follows immediately from Propositions \ref{prop:expectederror} and \ref{prop:asymptotics2}.

\begin{cor}\label{prop:mainresult}
As $m\to\infty$, we have
\[
\begin{split}
P_{n,m,r}(p)=&
\left[\frac12\erfc(\sqrt{m}\rho(r/m,p/(1-p)))\right.\\
&
+\left.\frac{A(r/m,p/(1-p))}{\sqrt{2m\pi}}e^{-m\rho(r/m,p/(1-p))^2}\right]^n(1+o(1))\;.
\end{split}
\]
\end{cor}

\begin{figure}[htb]
\vspace*{-2cm}
\begin{center}\includegraphics[angle=0,width=0.7\textwidth]{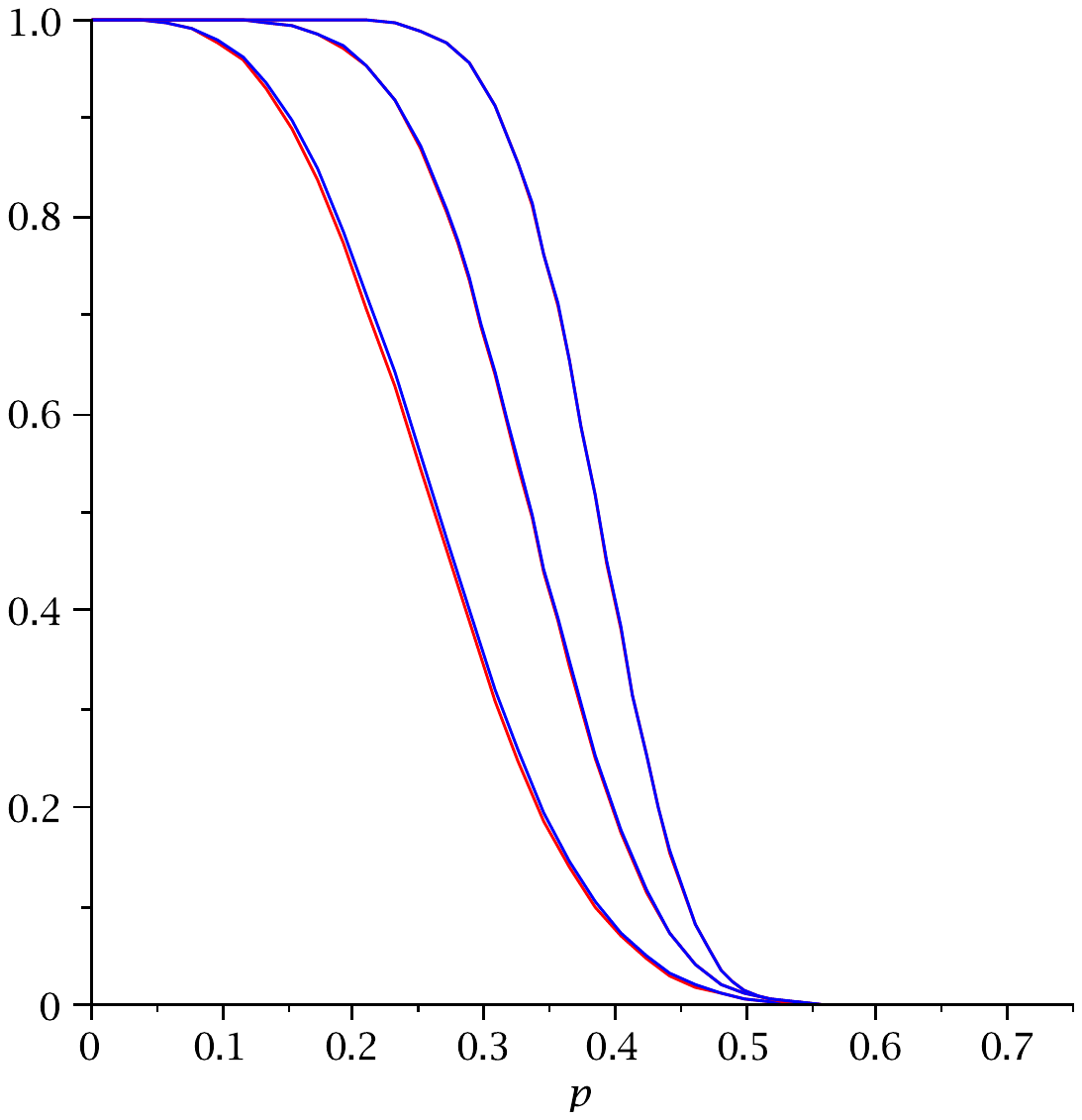}\end{center}
\vspace*{-2cm}
\caption{Shown are six curves: for $m=8$, $16$ and $32$ (from left to right), the respective curves for $P_{n,m,r}(p)$ and for the asymptotic result from Proposition \ref{prop:mainresult} are barely distinguishable.}
\label{figure3}
\end{figure}

\begin{example} \label{example:asy2example}
As we consider $r=\lfloor\frac{m-1}{2}\rfloor$, asymptotically $\beta=1/2$, but even for small values of $m$, such as $m=8$ and $\beta=3/8$, the expression in
Proposition \ref{prop:asymptotics2} provides a surprisingly accurate approximation, as can be seen from Figure \ref{figure3}.
\end{example}

\section{Conclusion} \label{section:conclusion}
In this note we have introduced a new algorithm for decoding the generalised hyperoctahedral group $C_m\wr S_n$.  If $n\ll m$, the performance of the algorithm is better from both the complexity perspective and from the number of correctible error patterns, when compared to the case $n\gg m$.

In particular, for large $m$ the complexity of the algorithm is $O(m\log m)$, whereas it is $O(n^2)$ for large $n$.  As is evident from Proposition \ref{prop:expectederror}, the number of correctible error patterns is a monotonically decreasing function of $n$ (this behaviour is demonstrated in Figure \ref{fig52}). On the other hand, using the properties of the complementary error function one can deduce from Corollary \ref{prop:mainresult} that for $p<r/m$ the number of correctible error patterns increases monotonically as a function of $m$ for $m$ sufficiently large (this behaviour is demonstrated in Figure~\ref{figure3}).

The results of the asymptotic analysis in Section \ref{section:asymptotics} give reasonable approximations even for moderately small values of $m$ and $n$. Extending Proposition \ref{prop:asymptotics2}, it is possible to give refined asymptotic estimates which are uniform in $m$ and $n$. 

As another direction, one could consider replacing the group $C_m \wr S_n$ with another wreath product, and modifying the algorithm to suit.  First, one could replace the symmetric group $S_n$ with another group $K$; however, this would give a much smaller number of codewords, and also the decoding algorithm would need to include a ``membership-testing'' algorithm (see Holt {\em et al.}\ \cite{Holt}) to check whether the decoded permutation was an element of $K$.  Second, one could replace the cyclic group with another group $H$; however, this would require a more sophisticated decoding process.

\section*{Acknowledgements}

The authors would like to thank P.~J.~Cameron and numerous others for reading the paper and providing helpful comments and suggestions.  The first author is a PIMS Postdoctoral Fellow at the University of Regina, and his research was supported in part by an EPSRC CASE studentship, sponsored by the UK Government Communications Headquarters (GCHQ).


\begin{thebibliography}{99}

\bibitem{btubb} R.~F.~Bailey, Uncoverings-by-bases for base-transitive permutation groups, {\em Des. Codes Cryptogr.} \textbf{41} (2006), 153--176.

\bibitem{ecpg} R.~F.~Bailey, Error-correcting codes from permutation groups, {\em Discrete Math.} {\bf 309} (2009), 4253--4265.

\bibitem{Blake74} I.~F.~Blake, Permutation codes for discrete channels, \emph{IEEE Trans. Inform. Theory} \textbf{20} (1974), 138--140.

\bibitem{BlakeCohenDeza79} I.~F.~Blake, G.~Cohen and M.~Deza, Coding with permutations, \emph{Information and Control} \textbf{43} (1979), 1--19.

\bibitem{Bleistein} N.~Bleistein and R.~A.~Handelsman, \emph{Asymptotic Expansions of Integrals}, Dover Publications, New York, 1986.

\bibitem{CameronPGbook} P.~J.~Cameron, \emph{Permutation Groups}, London Mathematical Society Student Texts (45), Cambridge University Press, Cambridge, 1999.

\bibitem{ChuColbournDukes04} W.~Chu, C.~J.~Colbourn and P.~Dukes, Constructions for permutation codes in powerline communications, \emph{Des. Codes Cryptogr.} \textbf{32} (2004), 51--64.

\bibitem{Flajolet} P.~Flajolet and R.~Sedgewick, {\em Analytic Combinatorics}, Cambridge University Press, Cambridge, 2009.

\bibitem{GAP} The GAP~Group, \emph{GAP -- Groups, Algorithms, and Programming, Version 4.4}; 2004, \url{http://www.gap-system.org}.

\bibitem{Holt} D.~F.~Holt, B.~Eick and E.~A.~O'Brien,  {\em Handbook of Computational Group Theory}, Chapman \& Hall, Boca Raton, 2005.

\bibitem{Huczynska06} S.~Huczynska, Powerline communications and the 36 officers problem, {\em Phil. Trans. Royal Soc. A} \textbf{364} (2006), 3199--3214.

\bibitem{TamoSchwartz10} I.~Tamo and M.~Schwartz, Correcting limited-magnitude errors in the rank-modulation scheme, {\em IEEE Trans. Inform. Theory} \textbf{56} (2010), 2551--2560.

\bibitem{Temme} N.~M.~Temme, Uniform asymptotic expansions of the incomplete gamma functions and the incomplete beta function, \emph{Math. Comp.} \textbf{29} (1975), 1109--1114.

\end{thebibliography}
\end{document}